\documentclass[12pt,a4paper,reqno]{amsart}
\usepackage{amsmath}
\usepackage{amsfonts}
\usepackage{amssymb}
\usepackage{latexsym}
\usepackage{graphicx}
\usepackage{color}
\usepackage{amscd}
\usepackage{epsfig}
\usepackage{cite}
\newtheorem{problem}{Problem}[section]
\newtheorem{definition}[problem]{Definition}
\newtheorem{lemma}[problem]{Lemma}
\newtheorem{theorem}[problem]{Theorem}

\newtheorem{corollary}[problem]{Corollary}

\title{$T$-adic exponential sums under diagonal base change}
\author{Chunlei Liu}\address{Department of Mathematical Sciences, Shanghai Jiao Tong
University, Shanghai 200240, P.R. China, E-mail: clliu@sjtu.edu.cn}
\begin{document}
\maketitle
\begin{abstract}
Twisted $T$-adic exponential sums are studied. The Hodge bound for
the $T$-adic Newton polygon of the $C$-function is established. As
an application, the behavior of the $L$-function under diagonal base
change is explicitly given.
\end{abstract}



\section{Introduction}
\subsection{Preliminaries} Let $\mathbb{F}_q$ be the field of
characteristic $p$ with $q$ elements, and
$\mathbb{Z}_q=W(\mathbb{F}_q)$. Let $T$ and $s$ be two independent
variables. In this subsection we are concerned with the ring
$\mathbb{Z}_q[[T]][[s]]$, elements of which are regarded as power
series in $s$ with coefficients in $\mathbb{Z}_q[[T]]$.

Let $\mathbb{Q}_p=\mathbb{Z}_p[\frac1p]$, $\overline{\mathbb{Q}}_p$
the algebraic closure of $\mathbb{Q}_p$, and
$\widehat{\overline{\mathbb{Q}}}_p$ the $p$-adic completion of
$\overline{\mathbb{Q}}_p$. \begin{definition}A (vertical)
specialization is a morphism $T\mapsto t$ from $\mathbb{Z}_q[[T]]$
into $\widehat{\overline{\mathbb{Q}}}_p$ with
$0\neq|t|_p<1$.\end{definition} We shall prove the vertical
specialization theorem.
\begin{theorem}[Vertical specialization]
Let $A(s,T)\in 1+s\mathbb{Z}_q[[T]][[s]]$ be a $T$-adic entrie
series in $s$. If $0\neq|t|_p<1$, then
$$t-adic\text{ NP of
}A(s,t)\geq T-adic\text{ NP of }A(s,T),$$where NP is the short for
Newton polygon. Moreover, the equality holds for one $t$ iff it
holds for all $t$.\end{theorem} By the vertical specialization, the
Newton polygon of a $T$-adic entire series in
$1+s\mathbb{Z}_q[[T]][[s]]$ goes up under vertical specialization,
and is stable under all specializations if it is stable under one
specialization.
\begin{definition}A $T$-adic entire series in $1+s\mathbb{Z}_q[[T]][[s]]$ is said to be stable
if its Newton polygon is stable under
specialization.\end{definition}
\begin{definition}[Tensor product]If
$$A(s,T)=\exp(-\sum\limits_{k=1}^{+\infty}a_k(T)\frac{s^k}{k}),$$
and
$$B(s,T)=\exp(-\sum\limits_{k=1}^{+\infty}b_k(T)\frac{s^k}{k}),$$
we define
$$A\otimes B(s,T)=\exp(-\sum\limits_{k=1}^{+\infty}a_k(T)b_k(T)\frac{s^k}{k}).$$
\end{definition}
We have the distribution law
$$(A_1A_2)\otimes B=(A_1\otimes B)(A\otimes B).$$
So, equipped with the usual multiplication and the new tensor
operation, the set of $T$-adic entire series in
$1+s\mathbb{Z}_q[[T]][[s]]$ becomes a ring. We shall prove that the
stable $T$-adic entire series form a subring.
\begin{theorem}The set of stable $T$-adic entire series
$1+s\mathbb{Z}_q[[T]][[s]]$ is closed under multiplication and
tensor operation.\end{theorem}

\subsection{Twisted $T$-adic exponential sums}In this subsection we
introduce $L$-functions of twisted $T$-adic exponential sums. The
theory of $T$-adic exponential sums without twists was developed by
Liu-Wan \cite{LW}.

Let $\mu_{q-1}$ be the group of $(q-1)$-th roots of unity in
$\mathbb{Z}_q$, $\omega:x\mapsto\hat{x}$ the Teichmuller character
of $\mathbb{F}_q^{\times}$ into $\mu_{q-1}$, $\chi=\omega^{-d}$ with
$s\in\mathbb{Z}^n/(q-1)$ a character of $(\mathbb{F}_q^{\times})^n$
into $\mu_{q-1}$, and $\chi_k=\chi\circ{\rm
Norm}_{\mathbb{F}_{q^k}/\mathbb{F}_q}$. Let $\psi(x)=(1+T)^x$ be the
quasi-character from $\mathbb{Z}_p$ to $\mathbb{Z}_p[[T]]^{\times}$,
and $\psi_q=\psi\circ\text{Tr}_{\mathbb{Z}_{q}/\mathbb{Z}_p}$. Let
$f\in\mu_{q-1}[x_1^{\pm1},\cdots,x_n^{\pm1}]$ be a non-constant
polynomial in $n$-variables with coefficients in $\mu_{q-1}$.
\begin{definition}The sum
$$S_{f}(T,\chi)=S_{f}(T,\chi,\mathbb{F}_q)=
\sum\limits_{x\in\mu_{q-1}^n}x^{-d} \psi_q\circ f(x)$$ is called a
twisted $T$-adic exponential sum. And the function
$$L_{f,\chi}(s,T)=L_{f,\chi}(s,T,\mathbb{F}_q)=
\exp(\sum\limits_{k=1}^{+\infty}S_{f}(T,\chi_k,\mathbb{F}_{q^k})\frac{s^k}{k})$$
is called an $L$-function of twisted exponential
sums.\end{definition} We have
$$L_{f,\chi}(s,T)=\prod_{x\in|\mathbb{G}_m^n\otimes\mathbb{F}_q|}
\frac{1}{1-\chi_{\deg(x)}(x)\psi_{\deg(x)}\circ
f(\hat{x})s^{\deg(x)}},$$ where $\mathbb{G}_m$ is the multiplicative
group $xy=1$.

That Euler product formula gives
$$L_{f,\chi}(s,T)\in 1+s\mathbb{Z}_q[[T]][[s]].$$

Define
$$C_{f,\chi}(s,T)=
\exp(\sum\limits_{k=1}^{+\infty}\frac{-1}{(q^k-1)^n}
S_{f}(T,\chi_k)\frac{s^k}{k}).$$ Call it a $C$-function of twisted
$T$-adic exponential sums. We have
$$L_{f,\chi}(s,T)=\prod_{i=0}^n
C_{f,\chi}(q^is,T)^{(-1)^{n-i+1}{n\choose i}} ,$$ and
$$C_{f,\chi}(s,T)=\prod_{j=0}^{+\infty}
L_{f,\chi}(q^js,T)^{(-1)^{n-1}{n+j-1\choose j}}.$$ So we have
$$C_{f,\chi}(s,T)\in1+s\mathbb{Z}_q[[T]][[s]].$$

We shall prove the analytic continuation of $C_{f,\chi}(s,T)$.
\begin{theorem}[Analytic continuation]The series $C_{f,\chi}(s,T)$ is $T$-adic entire
in $s$.
\end{theorem}

The analytic continuation of $C_{f,\chi}(s,T)$ immediately gives the
meromorphic continuation of $L_{f,\chi}(s,T)$.
\begin{theorem}[Meromorphic continuation]The series $L_{f,\chi}(s,T)$ is $T$-adic
meromorphic.
\end{theorem}
\begin{definition}The Laurent polynomial $f$ is said to be $\chi$-twisted stable
if $C_{f,\chi}(s,T)$ is stable $T$-adic entire series in
$s$.\end{definition}

Let $\zeta_{p^m}$ denote a primitive $p^m$-th root of unity. The
specialization
 $L_{f,\chi}(s,\zeta_{p^m}-1)$ is the $L$-function of twisted
 algebraic
exponential sums $S_{f}(\zeta_{p^m},\chi_k)$. These sums were
studied by Liu \cite{L}, with the $m=1$ case studied by
Adolphson-Sperber \cite{AS}.

 \begin{definition}Let $f(x)=\sum\limits_{u\in I}a_ux^u$ with
 $I\subseteq\mathbb{Z}$ and $a_u^{q-1}=1$. We define
$\triangle(f)$ to be the convex polytope
 in $\mathbb{R}^n$ generated by the origin and the vectors $u\in I$.\end{definition}
 \begin{definition}We call
$f$ non-degenerate if $\triangle(f)$ is of dimension $n$, and for
every closed face $\sigma\not\ni0$ of $\triangle(f)$, the system
$$\frac{\partial f_{\sigma}}{\partial x_1} \equiv\cdots \equiv\frac{\partial f_{\sigma}}{\partial
x_n} \equiv0(\mod p)$$ has no common zeros in
$(\overline{\mathbb{F}}_q^{\times})^n$, where
$f_{\sigma}=\sum\limits_{u\in\sigma}a_ux^u$. \end{definition}

Gelfand-Kapranov-Zelevinsky proved the following.
\begin{theorem}[\cite{GKZ}]Let
$\triangle\ni0$ be an integral convex polytope
 in $\mathbb{R}^n$. If $p$ is sufficiently large, and $f$ is a generic Laurent polynomial
 in $\triangle(f)=\triangle$, then $f$ is non-degenerate.\end{theorem}

By the above theorem, we are mainly concerned with non-degenerate
$f$. We have the following.
\begin{theorem}[\cite{L}]If $f$ is non-degenerate, then
$L_{f,\chi}(s,\zeta_{p^m}-1,\mathbb{F}_q)^{(-1)^{n-1}}$ is a
polynomial of degree $p^{n(m-1)}{\rm Vol}(\triangle(f))$.
\end{theorem}

For non-degenerate $f$, the determination of the Newton polygon of
$L_{f,\chi}(s,\zeta_{p^m}-1,\mathbb{F}_q)^{(-1)^{n-1}}$ is a
challenging problem. The case $m=1$ is already very difficult, let
alone the case $m>1$. However, from the vertical specialization
theorem, one can prove the following.
\begin{theorem}[Newton polygon for stable Laurent polynomials]
Suppose that $f$ is non-degenerate and $\chi$-twisted stable. Let
 $\lambda_1,\cdots,\lambda_r$ be the slopes of
 the $q$-adic Newton polygon of
$L_{f,\chi}(s,\zeta_p-1)^{(-1)^{n-1}}$. Then the slopes of the
$q$-adic Newton polygon of
$L_{f,\chi}(s,\zeta_{p^m}-1)^{(-1)^{n-1}}$ are the numbers
$$\frac{\lambda_i+j_1+j_2+\cdots+j_n}{p^{m-1}},$$
where $i=1,\cdots, r$, and each $j_k=0,1,\cdots,p^{m-1}-1$.
\end{theorem}
The non-degenerate condition for $f$ can be replaced by the
condition that the functions
$L_{f,\chi}(s,\zeta_{p}-1)^{(-1)^{n-1}}$ and
$L_{f,\chi}(s,\zeta_{p^m}-1)^{(-1)^{n-1}}$ are polynomials. The
above theorem reduces the determination of the Newton polygon of
$L_{f,\chi}(s,\zeta_{p^m}-1)$ to the $m=1$ case, provided that $f$
is non-degenerate and $\chi$-twisted stable. So, for non-degenerate
$f$, we are mainly concerned with the stability of $f$ and the
determination of the Newton polygon of
$L_{f,\chi}(s,\zeta_{p}-1)^{(-1)^{n-1}}$. We shall prove the
following stability criterion.
\begin{theorem}[Stability of ordinary Laurent polynomials]If $f$ is $\chi$-twisted ordinary, then it is
$\chi$-twisted stable.
 \end{theorem}

We now recall the notion of $\chi$-twisted ordinary Laurent
polynomial. Let $\triangle\ni0$ be an integral convex polytope in
$\mathbb{R}^n$, $C(\triangle)$ the cone generated by $\triangle$,
$M(\triangle)=C(\triangle)\cap\mathbb{Z}^n$, and $\deg_{\triangle}$
the degree function on $C(\triangle)$, which is $\mathbb{R}_+$
linear and takes the value $1$ on each face $\delta\not\ni0$. Let
$d\in\mathbb{Z}^n/(q-1)$, and
$$M_{d}(\triangle):=\frac{1}{q-1}(M(\triangle)\cap d).$$
\begin{definition}Let $b$ be the least positive integer such that $p^bd=d$. Order
elements of $\cup_{i=0}^{b-1}M_{p^id}(\triangle)$ so that
$$\deg_{\triangle}(x_1)\leq\deg_{\triangle}(x_2)\leq\cdots.$$
The infinite $d$-twisted Hodge polygon $H_{\triangle,d}^{\infty}$ of
$\triangle$ is the convex function on $\mathbb{R}_+$ with initial
value $0$ which is linear between consecutive integers and whose
slopes (between consecutive integers) are
$$\frac{\deg_{\triangle}(x_{bi+1})+\deg_{\triangle}(x_{bi+2})+\cdots+\deg_{\triangle}(x_{b(i+1)})}{b},\ i=0,1,\cdots.$$
\end{definition}
\begin{definition}If
$$T-adic \text{ NP of }
C_{f,\omega^{-d}}(s,T,\mathbb{F}_q)={\rm
ord}_p(q)(p-1)H_{\triangle(f),d}^{\infty},$$ then $f$ is called
$\omega^{-d}$-twisted ordinary.\end{definition}

\subsection{Exponential sums under diagonal base change}
In this subsection we introduce the exponential sums associated to
the tensor product of two Laurent polynomials.

\begin{definition}If
$g=\sum\limits_{v}b_vy^v\in\mu_{q-1}^n[y_1^{\pm1},\cdots,y_m^{\pm1}],$
we define $$f\otimes g=\sum\limits_{u,v}a_ub_v^uz^{u\otimes
v}\in\mu_{q-1}[z_{ij}^{\pm1},i=1,\cdots,n,j=1,\cdots,m],$$ and call
it a diagonal base change of $f$.\end{definition}

The number $q$ acts on the $m$-tuples $u=(u_1,\cdots,u_m)$ of
vectors in $\mathbb{Z}_{(p)}^n/\mathbb{Z}^n$ by multiplication. The
length of the orbit $u$ is denoted by $|u|$. The congruences
$u_1\otimes v_1+\cdots +u_m\otimes v_m\equiv \frac{d}{q-1}$, and
$(q^k-1)u_j\equiv0$ are defined on the orbit space
$q\setminus(\mathbb{Z}_{(p)}^n/\mathbb{Z}^n)^m$.

We shall prove the following.
\begin{theorem}Let $v_1,\cdots,v_m$ be an integral basis of $\mathbb{R}^m$,
and $g(y)=\sum\limits_jb_jy^{v_j}$ with $b_j\in\mu_{q-1}^n$. Then
$$C_{f\otimes g,\omega^{-d}}(s,T)=\prod_{\stackrel{(u_1,\cdots,u_m)\in
q\setminus(\mathbb{Z}_{(p)}^n/\mathbb{Z}^n)^m}{u_1\otimes v_1+\cdots
+u_m\otimes v_m\equiv \frac{d}{q-1}}}
\otimes_{j=1}^mC_{f,\omega^{-u_j(q^{|u|}-1)}}(s^{|u|}\prod_{j=1}^m
b_j^{-u_j(q^{|u|}-1)},T,\mathbb{F}_{q^r}).
$$\end{theorem}
As a function of $(u_1,\cdots,u_m)$, the tensor product on the
right-hand side of the equality is defined on the orbit space. And,
as the solutions of $u_1\otimes v_1+\cdots +u_m\otimes v_m\equiv0$,
under the map
$$(u_1,\cdots,u_m)\mapsto u_1\otimes v_1+\cdots +u_m\otimes v_m,$$
can be embedded into
$\mathbb{Z}^{n}\otimes\mathbb{Z}^m/\mathbb{Z}^n\otimes\sum\limits_{j=1}^m\mathbb{Z}v_j$,
the product on the right-hand side of the equality is a finite
product.

The above theorem has the following equivalent
form.\begin{theorem}Let $v_1,\cdots,v_m$ be an integral basis of
$\mathbb{R}^m$, and $g(y)=\sum\limits_jb_jy^{v_j}$ with
$b_j\in\mu_{q-1}^n$. Then

$$L_{f\otimes g,\omega^{-d}}(s,T)^{(-1)^{mn-1}}=\prod_{\stackrel{(u_1,\cdots,u_m)\in
q\setminus(\mathbb{Z}_{(p)}^n/\mathbb{Z}^n)^m}{u_1\otimes v_1+\cdots
+u_m\otimes v_m\equiv \frac{d}{q-1}}}
\otimes_{j=1}^mL_{f,\omega^{-u_j(q^{|u|}-1)}}(s^{|u|}\prod_{j=1}^m
b_j^{-u_j(q^{|u|}-1)},T,\mathbb{F}_{q^{|u|}})^{(-1)^{n-1}}.
$$\end{theorem}
By the above theorem, the Newton polygon of the $L$-function or
$C$-function of $T$-adic (resp. algebraic) exponential sums of
$f\otimes g$ is determined by that of $f$.

Combine the above with theorem the fact that the set of stable
$T$-adic entire series in $1+s\mathbb{Z}_q[[T]][[s]]$ is closed
under multiplication and tensor operation, we get the following.
\begin{corollary}Let $v_1,\cdots,v_m$
be an integral basis of $\mathbb{R}^m$, and
$g(x)=\sum\limits_{j=1}^mb_jx^{v_j}$ with $b_j\in\mu_{q-1}^n$. If
$f$ is $\chi$-twisted stable for all $\chi$, then so is $f\otimes
g$.\end{corollary}

\section{Analytic continuation}
In this section, we prove the analytic continuation of
$C_{f,\psi}(s,T,\mathbb{F}_q)$.

Define a new variable $\pi$ by the relation $E(\pi)=1+T$, where
$$E(\pi)=\exp(\sum_{i=0}^{\infty}\frac{\pi^{p^i}}{p^i}) \in 1+\pi{\mathbb Z}_p[[\pi]]$$
is the Artin-Hasse exponential series. Thus, $\pi$ is also a
$T$-adic uniformizer of ${\mathbb Q}_p((T))$.

Let $\triangle=\triangle(f)$, and $D$ the least common multiple of
the denominators of $\deg(\triangle)$. Write
$$L_d(\triangle)=\{\sum\limits_{u\in M_d(\triangle)}c_u\pi^{\deg(u)}x^u :
 c_u\in\mathbb{Z}_{q}[[\pi^{\frac{1}{D(q-1)}}]]\},$$
 and
$$B_d(\triangle)=\{\sum\limits_{u\in M_d(\triangle)}c_u\pi^{\deg(u)}x^u :
 c_u\in\mathbb{Z}_{q}[[\pi^{\frac{1}{D(q-1)}}]],\
 \text{ord}_T(c_u)\rightarrow+\infty
 \text{ if }\deg(u)\rightarrow+\infty\}.$$
Note that $L_d(\triangle)$ is stable under multiplication by
elements of $L_0(\triangle)$, and
 $$E_f(x) =\prod\limits_{a_u\neq0}E(\pi \hat{a}_u\hat{x}^u)\in L_0(\triangle).$$

Define
$$\phi:L_d(\Delta)\rightarrow L_{dp^{-1}}(\Delta),\
 \sum\limits_{u\in
M_d(\Delta)} c_ux^u\mapsto\sum\limits_{u\in M_{dp^{-1}}(\Delta)}
c_{pu}x^u.$$ Then the map $\phi\circ E_f$ sends $L_d$ to
$B_{dp^{-1}}$.

\begin{lemma}If $x^{p^a-1}=1$, then
$$E(\pi)^{x+x^p+\cdots+x^{p^{a-1}}}=E(\pi x)E(\pi x^p)\cdots E(\pi x^{p^{a-1}}).$$
\end{lemma} \proof  Since
$$\sum\limits_{j=0}^{a-1}x^{p^j}=\sum\limits_{j=0}^{a-1}x^{p^{j+i}},$$
we have
$$
E(\pi)^{x+x^p+\cdots+x^{p^{a-1}}}
=\exp(\sum_{i=0}^{\infty}\frac{\pi^{p^i}}{p^i}\sum\limits_{j=0}^{a-1}x^{p^{j+i}})=E(\pi
x)E(\pi x^p)\cdots E(\pi x^{p^{a-1}}).$$\endproof

The Galois group $\text{Gal}(\mathbb{Q}_q/\mathbb{Q}_p)$ is
generated by the Frobenius element $\sigma$, whose restriction to
$(q-1)$-th roots of unity is the $p$-power map. That Galois group
can act on $L(\Delta)$ by fixing $\pi^{\frac{1}{D(q-1)}}$ and
$x_1,\cdots,x_n$.
\begin{lemma}[Dwork's splitting lemma]If $q=p^a$, and
 $x\in(\mathbb{F}_{q^k}^{\times})^n$, then $$E(\pi)^{\text{Tr}_{\mathbb{Z}_{q^k}/\mathbb{Z}_p}
 (\sum\limits_u\hat{a}_u\hat{x}^u)}
 =\prod\limits_{i=0}^{ak-1}E_f^{\sigma^i}(\hat{x}^{p^i}).$$ \end{lemma} \proof We have
$$E(\pi)^{\text{Tr}_{\mathbb{Z}_{q^k}/\mathbb{Z}_p}
 (\sum\limits_u\hat{a}_u\hat{x}^u)}
 =\prod\limits_{a_u\neq0}E(\pi)^{\text{Tr}_{\mathbb{Z}_{q^k}/\mathbb{Z}_p}(\hat{a}_u\hat{x}^u)}$$
$$=\prod\limits_{a_u\neq0}\prod\limits_{i=0}^{ak-1}E(\pi(\hat{a}_u\hat{x}^u)^{p^i})=
\prod\limits_{i=0}^{ak-1}E_f^{\sigma^i}(\hat{x}^{p^i}).$$\endproof

Define $c(u,v)=\deg(u)+\deg(v)-\deg(u+v)\text{ if }u,v\in
C(\triangle)$. Then $c(u,v)\geq0$, and is zero if and only if $u$
and $v$ are cofacial. We call $c(u,v)$ the cofacial defect of $u$
and $v$.
\begin{lemma} Write
$$E_f(x) =\sum\limits_{u\in
M(\triangle)}\alpha_u(f)\pi^{\deg(u)}x^u.$$ Then, for $u\in
M_{d}(\Delta)$, we have
$$\phi\circ E_f(\pi^{\deg(u)}x^u)=\sum\limits_{w\in
M_{dp^{-1}}(\Delta)}\alpha_{pw-u}(f)\pi^{c(pw-u,u)}\pi^{(p-1)\deg(w)}\pi^{\deg(w)}x^w.$$
\end{lemma} \proof  Obvious.\endproof

Define $\phi_p:=\sigma^{-1}\circ\phi\circ E_f$, and
$\phi_{p^a}=\phi_p^{a}$. Then $\phi_{p^a}$ sends $B_d$ to
$B_{dp^{-a}}$, and
$$\phi_{p^{a}}=\sigma^{-a}\circ\phi^{a}\circ
\prod\limits_{i=0}^{a-1}E_{f}^{\sigma^i}(x^{p^i}).$$ It follows that
$\phi_q$ operates on $B_d$, and is linear over
$\mathbb{Z}_q[[\pi^{\frac{1}{D(q-1)}}]]$. Moreover,  by the last
lemma, it is completely continuous in the sense of \cite{Se}.

\begin{theorem}[Dwork's trace formula]\label{trace-formula}Suppose that $\chi=\omega^{-d}$. Then
$$S_{f,\chi}(T,\mathbb{F}_{q^k})
=(q^k-1)^n\text{Tr}_{B_d/\mathbb{Z}_{q}[[\pi^{\frac{1}{D(q-1)}}]]}(\phi_{q}^{k}),\
k=1,2\cdots.$$\end{theorem} \proof Suppose that $q=p^a$. Let
$g(x)\in B_d$. We have
$$\phi_q^{k}(g)=\phi^{ak}(g\prod\limits_{i=0}^{ak-1}E_f^{\sigma^i}(x^{p^i})).$$Write
$\prod\limits_{i=0}^{ak-1}E_f^{\sigma^i}(x^{p^i})=\sum\limits_{u\in
M(\Delta)}\beta_u x^u $. Then
$$\phi_q^{k}(\pi^{\deg(v)}x^v)=\sum\limits_{u\in
M_d(\Delta)}\beta_{q^ku-v}\pi^{\deg(v)}x^u.$$ So the trace of
$\phi_q^{k}$ on $B_d$ over
$\mathbb{Z}_{q}[[\pi^{\frac{1}{D(q-1)}}]]$ equals $\sum\limits_{u\in
M_d(\Delta)} \beta_{(q^k-1)u}$. But, by Dwork's splitting lemma, we
have
$$S_{f,\chi}(T,\mathbb{F}_{q^k})
=\sum\limits_{x_1^{q^k-1}=1,\cdots,x_n^{q^k-1}=1}
x^{-d(1+q+\cdots+q^{k-1})} \prod\limits_{i=0}^{ak-1}
E_f^{\sigma^i}(x^{p^i})=(q^k-1)^n\sum\limits_{u\in M_d(\Delta)}
\beta_{(q^k-1)u}.$$ The theorem now follows.\endproof

\begin{theorem}[Analytic trace formula]\label{analytic-trace-formula}
If $\chi=\omega^{-d}$, then
$$C_{f,\chi}(s,T,\mathbb{F}_q)={\rm det}_{\mathbb{Z}_{q}[[\pi^{\frac{1}{D(q-1)}}]]}(1-\phi_qs\mid
B_d).$$ In particular, $C_{f,\chi}(s,T,\mathbb{F}_q)$ is $T$-adic
analytic in $s$.
\end{theorem}

\proof This follows from the last theorem and the identity
$${\rm det}_{\mathbb{Z}_{q}[[\pi^{\frac{1}{D(q-1)}}]]}(1-\phi_qs\mid B_d)
=\exp(-\sum\limits_{k=1}^{+\infty}\text{Tr}_{B_d/\mathbb{Z}_{q}[[\pi^{\frac{1}{D(q-1)}}]]}(\phi_q^{k})\frac{s^k}{k}).$$
\endproof

\section{Hodge bound}
In this section, we prove the Hodge bound for the Newton polygon of
$C_{f,\chi}(s,T,\mathbb{F}_q)$. It will play an important role in
establishing the stability of ordinary Laurent polynomial.

 Let the Galois group
$\text{Gal}(\mathbb{Q}_q/\mathbb{Q}_p)$ act on
$\mathbb{Z}_q[[T]][[s]]$ by fixing $s$ and $T$.

\begin{lemma}We have
$$C_{f,\chi}(s,T,\mathbb{F}_q)^{\sigma}=C_{f,\chi^p}(s,T,\mathbb{F}_q).$$\end{lemma}
\proof Obvious.\endproof

\begin{corollary}Suppose that $q=p^a$ and $\chi=\omega^{-d}$.
Let $b$ be the least positive integer such that $p^bd=d$. Then, as
power series in $s$ with coefficients in $\mathbb{Z}_q[[T]]$,
$$\text{NP of }C_{f,\chi}(s,T,\mathbb{F}_q)^{ab}=
\text{NP of }{\rm
det}_{\mathbb{Z}_{p}[[\pi^{\frac{1}{D(q-1)}}]]}(1-\phi_qs\mid
\oplus_{i=0}^{b-1}B_{p^id}).$$
\end{corollary} \proof In fact, we
have
$${\rm det}_{\mathbb{Z}_p[[\pi^{\frac{1}{D(q-1)}}]]}(1-\phi_qs\mid
\oplus_{i=0}^{b-1}B_{p^id}) =\prod_{j=0}^{a-1}{\rm
det}_{\mathbb{Z}_q[[\pi^{\frac{1}{D(q-1)}}]]}(1-\phi_qs\mid
\oplus_{i=0}^{b-1}B_{p^id})^{\sigma^j}$$$$=\prod_{j=0}^{a-1}\prod_{i=0}^{b-1}{\rm
det}_{\mathbb{Z}_q[[\pi^{\frac{1}{D(q-1)}}]]}(1-\phi_qs\mid
B_{d})^{\sigma^{i+j}}.$$ The corollary now follows.\endproof

\begin{corollary}Suppose that $q=p^a$ and $\chi=\omega^{-d}$.
Let $b$ be the least positive integer such that $p^bd=d$. Then, as
power series in $s$ with coefficients in $\mathbb{Z}_q[[T]]$,
$$\text{NP of }C_{f,\chi}(s^a,T,\mathbb{F}_q)^{b}=
\text{NP of }{\rm
det}_{\mathbb{Z}_{p}[[\pi^{\frac{1}{D(q-1)}}]]}(1-\phi_ps\mid
\oplus_{i=0}^{b-1}B_{p^id}).$$
\end{corollary} \proof This follows from the identity
$${\rm det}_{\mathbb{Z}_p[[\pi^{\frac{1}{D(q-1)}}]]}(1-\phi_qs^a\mid
\oplus_{i=0}^{b-1}B_{p^id}) =\prod_{\zeta^a=1}{\rm
det}_{\mathbb{Z}_p[[\pi^{\frac{1}{D(q-1)}}]]}(1-\phi_p\zeta s\mid
\oplus_{i=0}^{b-1}B_{p^id}).$$\endproof

\begin{theorem}\label{hodge-for-psi}Suppose that $q=p^a$ and $b$ is the least positive integer such that $p^bd=d$. Then, as a
power series in $s$ with coefficients in $\mathbb{Z}_q[[T]]$, the
$T$-adic Newton polygon of ${\rm
det}_{\mathbb{Z}_{p}[[\pi^{\frac{1}{D(q-1)}}]]}(1-\phi_ps\mid
\oplus_{i=0}^{b-1}B_{p^id})$ lies above the convex polygon with
initial point $(0,0)$ and slopes $(p-1)\deg(w)$, where $w$ runs
through elements of $\cup_{i=0}^{b-1}M_{p^id}$ with multiplicity
$a$.\end{theorem}

\proof Choose $\zeta\in\mathbb{Z}_q^{\times}$ such that
$\zeta^{\sigma^i}$, $i=0,1,\cdots, a-1$ be a basis of $\mathbb{Z}_q$
over $\mathbb{Z}_p$. Write
$$\alpha_u(f)=\sum\limits_{i=0}^{a-1}\alpha_{u,i}(f)\zeta^{\sigma^i},\
\alpha_{u,i}(f)\in\mathbb{Z}_p[[\pi^{1/D}]].$$ Then, for $u\in
M_{p^id}(\triangle)$, we have
$$\phi_p(\zeta^{\sigma^j}\pi^{\deg(u)}x^u)=\sum\limits_{i=0}^{a-1}\sum\limits_{w\in
M_{dp^{i-1}}(\Delta)}
\alpha_{pw-u,i-j+1}(f)\pi^{c(pw-u,u)}\pi^{(p-1)\deg(w)}\zeta^{\sigma^i}\pi^{\deg(w)}x^w.$$
So, the matrix of $\phi_p$ over
$\mathbb{Z}_p[[\pi^{\frac{1}{D(q-1)}}]]$ with respect to the basis
$\{\zeta^{q^j}\pi^{\deg(u)}x^u\}_{0\leq j<a,u\in
\oplus_{i=0}^{b-1}B_{p^id}(\Delta)}$ is
$$A=(\alpha_{pw-u,i-j+1}(f)\pi^{c(pw-u,u)}\pi^{(p-1)\deg(w)})_{(i,w),(j,u)}.$$
It follows that, the $T$-adic Newton polygon of ${\rm
det}_{\mathbb{Z}_{p}[[\pi^{\frac{1}{D(q-1)}}]]}(1-\phi_ps\mid
\oplus_{i=0}^{b-1}B_{p^id})$ lies above the convex polygon with
initial point $(0,0)$ and slopes $(p-1)\deg(w)$, where $w$ runs
through elements of $\cup_{i=0}^{b-1}M_{p^id}$ with multiplicity
$a$.\endproof

\begin{corollary}\label{hodge2}Suppose that $q=p^a$ and $\chi=\omega^{-d}$.
Let $b$ be the least positive integer such that $p^bd=d$. Then, as a
power series in $s$ with coefficients in $\mathbb{Z}_q[[T]]$, the
$T$-adic Newton polygon of $C_{f,\chi}(s,T,\mathbb{F}_q)^{b}$ lies
above the convex polygon with initial point $(0,0)$ and slopes
$$a(p-1)\deg(w),\ w\in\cup_{i=0}^{b-1}M_{p^id}.$$\end{corollary}

\proof Obvious.\endproof

\begin{theorem}[Hodge bound]\label{hodge-bound}Suppose that $\chi=\omega^{-d}$, and
$q=p^a$. Then $$T-adic \text{ NP of } C_{f,\chi}(s,T)\geq
a(p-1)H_{\triangle(f),d}^{\infty},$$ where NP is the short for
Newton polygon,
 and $H_{\triangle(f),d}^{\infty}$ is the infinite $d$-twisted Hodge polygon
of $\triangle(f)$.\end{theorem}

\proof Obvious.\endproof

 \begin{definition}If $\chi=\omega^{-d}$,
$q=p^a$, $0\neq|t|_p<1$, and $$t-adic\text{ NP of
}C_{f,\chi}(s,t)=a(p-1)H_{\triangle(f),d}^{\infty},$$ then $f$ is
said to be $\chi$-twisted ordinary.\end{definition}
\section{Vertical specialization and stability} In this section we prove
the vertical specialization theorem, the theorem for the Newton
polygon of stable Laurent polynomials, and the stability of ordinary
Laurent polynomials.
\begin{theorem}[Vertical specialization]
Let $A(s,T)\in 1+s\mathbb{Z}_q[[T]][[s]]$ be a $T$-adic entrie
series in $s$. If $0\neq|t|_p<1$, then
$$t-adic\text{ NP of
}A(s,t)\geq T-adic\text{ NP of }A(s,T),$$where NP is the short for
Newton polygon. Moreover, the equality holds for one $t$ iff it
holds for all $t$. \end{theorem} \proof Write
 $$A(s,T)=\sum\limits_{i=0}^{\infty}a_i(T)s^i.$$
 The inequality follows from the fact that
$a_i(T)\in\mathbb{Z}_q[[T]]$. Moreover, $$t-adic\text{ NP of
}A(s,t)=T-adic\text{ NP of }A(s,T)$$if and only if
$$a_i(T)\in T^e\mathbb{Z}_q[[T]]^{\times}$$ for every turning point
$(i,e)$ of the $T$-adic Newton polygon of $A(s,T)$.  It follows that
the equality holds for one $t$ iff it holds for all $t$.\endproof
\begin{theorem}[Newton polygon for stable Laurent polynomials]
Suppose that $f$ is non-degenrate and $\chi$-twisted stable. Let
 $\lambda_1,\cdots,\lambda_r$ be the slopes of
 the $q$-adic Newton polygon of $L_{f,\chi}(s,\zeta_p-1)^{(-1)^{n-1}}$.
 Then the $q$-adic orders of the reciprocal roots of
$L_{f,\chi}(s,\zeta_{p^m}-1)^{(-1)^{n-1}}$ are the numbers
$$\frac{\lambda_i+j_1+j_2+\cdots+j_n}{p^{m-1}},$$
where $i=1,\cdots, r$, and each $j_k=0,1,\cdots,p^{m-1}-1$.
\end{theorem}  \proof Apply the relationship between
the $L$-function and the $C$-function, we see that the $q$-adic
orders of the reciprocal roots of $C_{f,\chi}(s,\zeta_{p}-1)$ are
the numbers $$\lambda_i+j_1+j_2+\cdots+j_n,$$ where $i=1,\cdots, r$,
and each $j_k=0,1,\cdots$. So the $(\zeta_p-1)$-adic orders of the
reciprocal roots of $C_{f,\chi}(s,\zeta_{p}-1)$ are the numbers
$${\rm ord}_{\zeta_p-1}(q)(\lambda_i+j_1+j_2+\cdots+j_n),$$ where $i=1,\cdots, r$, and each $j_k=0,1,\cdots$. Apply the
$\chi$-twisted stability of $f$, we see that the
$(\zeta_{p^m}-1)$-adic orders of the reciprocal roots of
$C_{f,\chi}(s,\zeta_{p^m}-1)$ are the numbers
$${\rm ord}_{\zeta_p-1}(q)(\lambda_i+j_1+j_2+\cdots+j_n),$$ where $i=1,\cdots, r$, and each $j_k=0,1,\cdots$. So the $q$-adic orders
of the reciprocal roots of $C_{f,\chi}(s,\zeta_{p^m}-1)$ are the
numbers
$$\frac{\lambda_i+j_1+j_2+\cdots+j_n}{p^{m-1}},$$ where $i=1,\cdots,r$, and each $j_k=0,1,\cdots$. Apply the relationship
between the $C$-function and the $L$-function, we see that the
$q$-adic orders of the reciprocal roots of
$L_{f,\chi}(s,\zeta_{p^m}-1)$ are the numbers
$$\frac{\lambda_i+j_1+j_2+\cdots+j_n}{p^{m-1}},$$ where $i=1,\cdots,r$, and each $j_k=0,1,\cdots,p^{m-1}-1$.\endproof

\begin{theorem}[Specialization of the Hodge bound]If $\chi=\omega^{-d}$,
$q=p^a$, and $0\neq|t|_p<1$, then $$t-adic\text{ NP of
}C_{f,\chi}(s,t)\geq T-adic\text{ NP of }C_{f,\chi}(s,T)\geq
a(p-1)H_{\triangle(f),d}^{\infty}.$$Moreover, the equalities hold
for one $t$ iff they hold for all $t$.\end{theorem} \proof Just
combine the Hodge bound for the Newton polygon of $C_{f,\chi}(s,T)$
with the vertical specialization theorem. \endproof
\begin{theorem}[Stability of ordinary Laurent polynomials]
If $f$ is $\chi$-twisted ordinary, then it is $\chi$-twisted stable,
and $\chi$-twisted $T$-adic ordinary.
 \end{theorem}
 \proof Obvious.\endproof

\begin{definition}Let $\alpha_1,\alpha_2,\cdots$ be the slopes
of the infinite $d$-twisted Hodge polygon of $\triangle$. Then
$$(1-t)^n\sum\limits_{i}t^{\alpha_i}=\sum\limits_{i=1}^{n!\text{Vol}(\triangle)}t^{w_i}.$$
The $d$-twisted Hodge polygon $H_{\triangle,d}$ of $\triangle$ is
the convex function on $[0,n!\text{Vol}(\triangle)]$ with initial
value $0$ which is linear between consecutive integers and whose
slopes (between consecutive integers) are $w_i$,
$i=1,\cdots,n!\text{Vol}(\triangle)$.
\end{definition}
\begin{theorem}[Newton polygon for ordinary Laurent polynomials]
Let $f$ be non-degenerate, $\chi=\omega^{-d}$ and $m\geq1$. Then $f$
is $\chi$-twisted ordinary if and only if
$$q-adic \text{ NP of }
L_{f,\chi}(s,\zeta_{p^m}-1)^{(-1)^{n-1}}=H_{p^{m-1}\triangle(f),d}.$$\end{theorem}
The non-degenerate condition for $f$ can be replaced by the
condition that the function
$L_{f,\chi}(s,\zeta_{p^m}-1)^{(-1)^{n-1}}$ is a polynomial. \proof
In fact, $f$ is $\chi$-twisted ordinary if and only if
$$(\zeta_p-1)-adic\text{ NP of }C_{f,\chi}(s,\zeta_p-1)={\rm ord}_p(q)(p-1)H_{\triangle(f),d}^{\infty},$$
if and only if $$(\zeta_{p^m}-1)-adic\text{ NP of
}C_{f,\chi}(s,\zeta_{p^m}-1)={\rm
ord}_p(q)(p-1)H_{\triangle(f),d}^{\infty},$$ if and only if
$$q-adic\text{ NP of }C_{f,\chi}(s,\zeta_{p^m}-1)
=\frac{1}{p^{m-1}}H_{\triangle(f),d}^{\infty}
=H_{p^{m-1}\triangle(f),d}^{\infty},$$ if and only if
$$q-adic \text{
NP of }
L_{f,\chi}(s,\zeta_{p^m}-1)^{(-1)^{n-1}}=H_{p^{m-1}\triangle(f),d}.$$.\endproof

\begin{lemma}[Hasse-Davenport relation]
$$L_{x,\chi}(s,\zeta_p,\mathbb{F}_{q})=1+G(\zeta_p,\chi,\mathbb{F}_q)s.$$\end{lemma}
\proof This follows from the the following classical
formulation:$$G(\zeta_p,\chi_k,\mathbb{F}_{q^k})=(-1)^{k-1}G(\zeta_p,\chi,\mathbb{F}_{q})^k.$$\endproof
\begin{definition}Let $q=p^a$, and $d\in\mathbb{Z}/(q-1)$. We define
$$\sigma_q(d)=(p-1)\sum\limits_{i=0}^{a-1}\{\frac{p^id}{q-1}\}.$$
\end{definition}
\begin{theorem}[Stickelberger theorem for Gauss sums]
The polynomial $f(x)=x$ is $\chi$-twisted ordinary for all
$\chi$.\end{theorem}

\proof We have $\triangle=\triangle(f)=[0,1]$,
$C(\triangle)=\mathbb{R}_+$, and $\deg_{\triangle}(u)=u$. Let
$q=p^a$, and $\chi=\omega^{-d}$. We have $M(\triangle)=\mathbb{N}$,
and $M_{p^id}(\triangle)=\{\frac{p^id}{q-1}\}+\mathbb{N}$. It
follows that the infinite $d$-twisted Hodge polygon of $\triangle$
has slopes
$$\frac{\sigma_q(d)}{a(p-1)}+k,\ k=0,1,\cdots.$$ So the
finite $d$-twisted Hodge polygon of $\triangle$ has only one slope
$\frac{\sigma_q(d)}{a(p-1)}$. By the Hasse-Davenport relation and
the classical Stickelberger theorem, the Newton polygon of the
$L$-function $L_{x,\chi}(s,\zeta_p-1)$ also has only one slope
$\frac{\sigma_q(d)}{a(p-1)}$. Therefore the polynomial $f(x)=x$ is
$\chi$-twisted ordinary.\endproof

\begin{theorem}[Stickelberger theorem for Gauss-Heilbronn sums]
Let $q=p^a$, and $\chi=\omega^{-d}$. Then $q$-adic orders of the
reciprocal zeros of the $L$-function $L_{x,\chi}(s,\zeta_{p^m}-1)$
of the Gauss-Heilbronn sums
$G(\zeta_{p^m}-1,\chi_k,\mathbb{F}_{q^k})$ are
$$\frac{\sigma_q(d)}{a(p-1)p^{m-1}}+\frac{k}{p^{m-1}},\ k=0,1,\cdots,p^{m-1}-1.$$\end{theorem}
The above theorem was proved by Blache \cite{B}, and Liu \cite{L}.
But the proof here is much simpler.

\proof Apply the Stickelberger theorem for Gauss sums and the
theorem on the Newton polygon for ordinary $p$-power order
exponential sums, we get $$q-adic \text{ NP of }
L_{x,\chi}(s,\zeta_{p^m}-1)=H_{p^{m-1}\triangle,d}$$ with
$\triangle=[0,1]$. In the proof of the Stickelberger theorem for
Gauss sums, we show that
$$M_{p^id}(\triangle)=\{\frac{p^id}{q-1}\}+\mathbb{N}.$$
It follows that the infinite $d$-twisted Hodge polygon of
$p^{m-1}\triangle$ has slopes
$$\frac{\sigma_q(d)}{a(p-1)p^{m-1}}+\frac{k}{p^{m-1}},\ k=0,1,\cdots.$$ So the
finite $d$-twisted Hodge polygon of $p^{m-1}\triangle$ has slopes
$$\frac{\sigma_q(d)}{a(p-1)p^{m-1}}+\frac{k}{p^{m-1}},\
k=0,1,\cdots,p^{m-1}-1.$$ The theorem now follows.
\endproof
\section{Stable $T$-adic entire series}In this section
we prove that the set of stable $T$-adic entire series in
$1+s\mathbb{Z}_q[[T]][[s]]$ is closed under multiplication and
tensor operation.
\begin{lemma}[Weierstrass preparation theorem]Let $A(s,T)\in\mathbb{Z}_q[[T]]\langle s\rangle$
be a $T$-adically strictly convergent power series in $s$ with
unitary constant term. Suppose that $A(s,T)(\mod
T)\in\mathbb{Z}_q[s]$ is a unitary polynomial of degree $n$. Then
$$A(s,T)=u(s,T)B(s,T),$$ where $u(s,T)\in1+T\mathbb{Z}_q[[T]]\langle
s\rangle$, and $B(s,T)\in\mathbb{Z}_q[[T]][s]$ is a monic polynomial
of degree $n$ with unitary constant term.
\end{lemma}
\proof Let $A_0(s)=A(s,T)(\mod T)$, and $\alpha={\rm
ord}_T(A(s,T)-A_0(s))$. Then $$\mathbb{Z}_q[[T]]\langle
s\rangle/(T^{\alpha},A(s,T))=\mathbb{Z}_q[[T]]\langle
s\rangle/(T^{\alpha},A_0(s)),$$ and is generated as
$\mathbb{Z}_q[[T]]$-module by $1,s,\cdots,s^{n-1}$. In particular.
$$v_0=s^n=\sum\limits_{j=0}^{n-1}a_{1j}s^j+v_1T^{\alpha}+w_1A(s,T),\
v_1\in\mathbb{Z}_q[[T]]\langle
s\rangle,w_1\in1+T\mathbb{Z}_q[[T]]\langle s\rangle.$$ By induction,
we can construct sequences $v_i,w_i\in\mathbb{Z}_q[[T]]\langle
s\rangle$ so
that$$v_{i-1}=\sum\limits_{j=0}^{n-1}a_{ij}s^j+v_iT^{\alpha}+w_iA(s
,T).$$ We have
$$\sum\limits_{i=1}^{\infty}v_{i-1}T^{(i-1)\alpha}
=\sum\limits_{j=0}^{n-1}s^j\sum\limits_{i=1}^{\infty}a_{ij}T^{(i-1)\alpha}
+\sum\limits_{i=1}^{\infty}v_iT^{i\alpha}+A(s
,T)\sum\limits_{i=1}^{\infty}w_iT^{(i-1)\alpha}.$$ So
$$s^n-\sum\limits_{j=0}^{n-1}s^j\sum\limits_{i=1}^{\infty}a_{ij}T^{(i-1)\alpha}=A(s
,T)\sum\limits_{i=1}^{\infty}w_iT^{(i-1)\alpha}.$$ Since
$w=\sum\limits_{i=1}^{\infty}w_iT^{(i-1)\alpha}\in1+T\mathbb{Z}_q[[T]]\langle
s\rangle$, we have $u(s,T)=w^{-1}\in1+T\mathbb{Z}_q[[T]]\langle
s\rangle$. Set
$B(s,T)=s^n-\sum\limits_{j=0}^{n-1}s^j\sum\limits_{i=1}^{\infty}a_{ij}T^{(i-1)\alpha}\in\mathbb{Z}_q[[T]][s]$,
we get $$A(s,T)=u(s,T)B(s,T).$$
\endproof

\begin{theorem}[Weierstrass factorization theorem]Let $A(s,T)\in1+s\mathbb{Z}_q[[T]]\langle s\rangle$
be a $T$-adically entire power series in $s$, whose Newton polygon
has slopes $\lambda_i$ of horizontal length $n_i$. Then
$$A(s,T)=\prod_{i=1}^{\infty}A_i(s),$$
where
$A_i(s)\in\mathbb{Z}_q[[T^{\lambda_1},\cdots,T^{\lambda_i}]][s]$ is
polynomial of degree $n_i$ with unitary constant term and linear
Newton polygon.
\end{theorem}
\proof Apply the Weierstrass preparation theorem to construct
$A_i(s)$ inductively.\endproof
\begin{lemma}Let $A(s,T)$ and $B(s,T)$ be two $T$-adic entire power series in
$1+s\mathbb{Z}_q[[T]][[s]]$. Suppose that
$$A(s,T)=\prod_{\alpha\in I}(1-\alpha s),$$
and $$B(s,T)=\prod_{\beta\in J}(1-\beta s).$$ Then $$A\otimes
B(s,T)=\prod_{\alpha\in I,\beta\in J}(1-\alpha\beta s).$$
\end{lemma}
\proof We have
$$A(s,T)=\exp(-\sum\limits_{k=1}^{+\infty}\frac{s^k}{k}
\sum\limits_{\alpha\in I}\alpha^k),$$ and
$$B(s,T)=\exp(-\sum\limits_{k=1}^{+\infty}\frac{s^k}{k}
\sum\limits_{\beta\in J}\beta^k).$$ So
$$A\otimes B(s,T)=\exp(-\sum\limits_{k=1}^{+\infty}\frac{s^k}{k}
\sum\limits_{\alpha\in I,\beta\in
J}\alpha^k\beta^k)=\prod_{\alpha\in I,\beta\in J}(1-\alpha\beta
s).$$\endproof
\begin{lemma}Let $A(s,T)$ and $B(s,T)$ be two $T$-adic entire power series in
$1+s\mathbb{Z}_q[[T]][[s]]$. Then $A(s,T)B(s,T)$ is stable iff both
$A(s,T)$ and $B(s,T)$ are stable.
\end{lemma} \proof Write
$$A(s,T)=\sum\limits_nA_n(T)s^n,$$
$$B(s,T)=\sum\limits_nB_n(T)s^n,$$
and $$A(s,T)B(s,T)=\sum\limits_nC_n(T)s^n.$$Let $\{\alpha\}$ be the
set of the reciprocal zeros of $A(s,T)$, and $\{\beta\}$ the set of
reciprocal zeros of $B(s,T)$. Let $(n,e)$ be a turning points of the
Newton polygon of $A(s,T)B(s,T)$. Then $n=\sum\limits_{{\rm
ord}_T(\alpha)\leq r}1+\sum\limits_{{\rm ord}_T(\beta)\leq r}1$ for
some $r\in\mathbb{R}_+$, and $$C_n(T)\equiv (-1)^n\prod_{{\rm
ord}_T(\alpha)\leq r}\alpha\prod_{{\rm ord}_T(\beta)\leq
r}\beta\equiv A_{n_1}(T)B_{n_2}(T)(\mod T^{>e}),$$ where
$n_1=\sum\limits_{{\rm ord}_T(\alpha)\leq r}1$, and
$n_2=\sum\limits_{{\rm ord}_T(\beta)\leq r}1$. Note that the
$T$-adic order of $C_n(T)$ is stable under specialization iff both
the $T$-adic orders of $A_{n_1}(T)$ and $B_{n_2}(T)$ are stable
under specialization. The lemma now follows.
\endproof
\begin{theorem}The set of stable $T$-adic entire series in
$1+s\mathbb{Z}_q[[T]][[s]]$ is closed under tensor
operation.\end{theorem} \proof Let $A(s,T)$ and $B(s,T)$ be two
stable $T$-adic entire power series in $1+s\mathbb{Z}_q[[T]][[s]]$.
By the Weierstrass factorization theorem, and the last lemma, we may
assume that $A$ and $B$ are polynomials with linear Newton polygon.
Let $\{\alpha_1,\cdots,\alpha_m\}$ be the set of the reciprocal
zeros of $A(s,T)$, and $\{\beta_1,\cdots,\beta_n\}$ the set of
reciprocal zeros of $B(s,T)$. Then the leading term of $A\otimes B$
is
$$C_{mn}(T)=\prod_{i=1}^m\prod_{j=1}^n(-\alpha_i\beta_j)=
(-1)^{mn}A_m(T)^nB_{n}^m(T),$$ where $A_m$ is the leading
coefficient of $A(s,T)$, and $B_n$ is the leading coefficient of
$B(s,T)$. Since the $T$-adic orders of $A_{m}(T)$ and $B_{n}(T)$ do
not go up under specialization, so does the $T$-adic order of
$C_{mn}(T)$. The theorem is proved.
\endproof
\section{Exponential sums under the tensor operation}
In this section we explore Wan's method \cite{W}, and study the
exponential sums associated to the tensor product of two Laurent
polynomials.
\begin{lemma}Let $v_1,\cdots,v_m$ be an integral basis of
$\mathbb{R}^m$, and $g(y)=\sum\limits_jb_jy^{v_j}$ with
$b_j\in\mu_{q-1}^n$. Then

$$S_{f\otimes g}(T,\omega^{-d})
=\sum\limits_{\stackrel{u_1,\cdots,u_m\in\mathbb{Z}^n/(q-1)}{u_1\otimes
v_1+\cdots +u_m\otimes v_m\equiv d}}
\prod_{j=1}^mb_j^{u_j}\prod_{j=1}^m
S_f(T,\omega^{-u_j},\mathbb{F}_q)
$$\end{lemma}
\proof Since
$$S_{f}(T,\omega^{-u},\mathbb{F}_q)
=\sum\limits_{\alpha\in\mu_{q-1}^n}\alpha^{-u}\psi_q\circ
f(\alpha),$$ we have
$$\psi_q\circ f(\alpha)=\frac{1}{(q-1)^n}\sum\limits_{u\in\mathbb{Z}^n/(q-1)}
\alpha^uS_f(T,\omega^{-u},\mathbb{F}_q).$$ So
$$\psi_q(f\otimes b_jy^{v_j}))
=\frac{1}{(q-1)^n}\sum\limits_{u\in\mathbb{Z}^n/(q-1)}
b_j^uz^{u\otimes v_j}S_f(T,\omega^{-u},\mathbb{F}_q).$$ Thus,
$$\psi_q(f\otimes g(z))
=\frac{1}{(q-1)^{mn}}\prod_j\sum\limits_{u\in\mathbb{Z}^n/(q-1)}
b_j^uz^{u\otimes v_j}S_f(T,\omega^{-u},\mathbb{F}_q).$$ Therfore
$$S_{f\otimes g}(T,\omega^{-d})=\sum\limits_{z\in\mu_{q-1}^{mn}}
\frac{z^{-d}}{(q-1)^{mn}}\sum\limits_{u_1,\cdots,u_m\in\mathbb{Z}^n/(q-1)}
z^{u_1\otimes v_1+\cdots+u_m\otimes
v_m}\prod_{j=1}^mb_j^{u_j}S_f(T,\omega^{-u_j}).$$ Change the order
of summation, we get the desired formula.\endproof

\begin{corollary}Let $v_1,\cdots,v_m$ be an integral basis of $\mathbb{R}^m$,
and $g(y)=\sum\limits_jb_jy^{v_j}$ with $b_j\in\mu_{q-1}^n$. Then

$$S_{f\otimes g}(T,\omega^{-\frac{d}{q-1}(q^k-1)},\mathbb{F}_{q^k})
=\sum\limits_{\stackrel{u_1,\cdots,u_m\in\mathbb{Z}_{(p)}^n/\mathbb{Z}^n}{u_1\otimes
v_1+\cdots +u_m\otimes v_m\equiv \frac{d}{q-1},(q^k-1)u_j\equiv0}}
\prod_{j=1}^mb_j^{u_j(q^k-1)}\prod_{j=1}^m
S_f(T,\omega^{-u_j(q^k-1)},\mathbb{F}_{q^k}).
$$\end{corollary}
\proof Just scale the variables $u_j$ in the last lemma. \endproof

\begin{lemma}We have
$$S_f(T,\chi^q,\mathbb{F}_{q^k})=S_f(T,\chi,\mathbb{F}_{q^k}).$$\end{lemma}
\proof Since $\sigma:x\mapsto x^q$ is an automorphism of
$\mu_{q^k-1}$, and extends to be an element of ${\rm
Gal}(\mathbb{Q}_{q^k}/\mathbb{Q}_q)$, we have
$$S_f(T,\chi^q,\mathbb{F}_{q^k})
=\sum\limits_{x\in\mu_{q^k-1}^n}\chi(x^q)\psi_{q^k}\circ f(x)
=\sum\limits_{x\in\mu_{q^k-1}^n}\chi(x)\psi_{q^k}\circ
f(x^{\sigma^{-1}}).$$ Note that
$$\psi_{q^k}\circ f(x^{\sigma^{-1}})
=\psi_q\circ{\rm
Tr}_{\mathbb{Q}_{q^k}/\mathbb{Q}_q}(f(x)^{\sigma^{-1}})=\psi_q\circ{\rm
Tr}_{\mathbb{Q}_{q^k}/\mathbb{Q}_q}(f(x))=\psi_{q^k}\circ f(x).$$The
lemma now follows.\endproof

By the above lemma, the product
$$\prod_{j=1}^mb_j^{-u_j(q^k-1)}\prod_{j=1}^m
S_f(T,\omega^{-u_j(q^k-1)},\mathbb{F}_{q^k})),$$ as a function of
$(u_1,\cdots,u_m)$, is also defined on the orbit space
$q\setminus(\mathbb{Z}_{(p)}^n/\mathbb{Z}^n)^m$. So we can restate
the last corollary as follows.

\begin{corollary}Let $v_1,\cdots,v_m$ be an integral basis of $\mathbb{R}^m$,
and $g(y)=\sum\limits_jb_jy^{v_j}$ with $b_j\in\mu_{q-1}^n$. Then

$$S_{f\otimes g}(T,\omega^{-\frac{d}{q-1}(q^k-1)},\mathbb{F}_{q^k})
=\sum\limits_{\stackrel{(u_1,\cdots,u_m)\in
q\setminus(\mathbb{Z}_{(p)}^n/\mathbb{Z}^n)^m}{u_1\otimes v_1+\cdots
+u_m\otimes v_m\equiv \frac{d}{q-1},(q^k-1)u_j\equiv0}}
|u|\prod_{j=1}^mb_j^{u_j(q^k-1)}\prod_{j=1}^m
S_f(T,\omega^{-u_j(q^k-1)},\mathbb{F}_{q^k}).
$$\end{corollary}

\begin{theorem}Let $v_1,\cdots,v_m$ be an integral basis of $\mathbb{R}^m$,
and $g(y)=\sum\limits_jb_jy^{v_j}$ with $b_j\in\mu_{q-1}^n$. Then

$$C_{f\otimes g,\omega^{-d}}(s,T)=\prod_{\stackrel{(u_1,\cdots,u_m)\in
q\setminus(\mathbb{Z}_{(p)}^n/\mathbb{Z}^n)^m}{u_1\otimes v_1+\cdots
+u_m\otimes v_m\equiv \frac{d}{q-1}}}
\otimes_{j=1}^mC_{f,\omega^{-u_j(q^{|u|}-1)}}(s^{|u|}\prod_{j=1}^m
b_j^{-u_j(q^{|u|}-1)},T,\mathbb{F}_{q^{|u|}}).
$$\end{theorem}
\proof We have
$$\sum\limits_{k=1}^{\infty}\frac{-1}{(q^k-1)^{mn}}S_{f\otimes g}(T,\omega^{-\frac{d}{q-1}(q^k-1)},\mathbb{F}_{q^k})\frac{s^k}{k}
$$$$=\sum\limits_{k=1}^{\infty}\frac{-1}{(q^k-1)^{mn}}\frac{s^k}{k}\sum\limits_{\stackrel{(u_1,\cdots,u_m)\in
q\setminus(\mathbb{Z}_{(p)}^n/\mathbb{Z}^n)^m}{u_1\otimes v_1+\cdots
+u_m\otimes v_m\equiv \frac{d}{q-1},(q^k-1)u_j\equiv0}}
|u|\prod_{j=1}^mb_j^{u_j(q^k-1)}\prod_{j=1}^m
S_f(T,\omega^{-u_j(q^k-1)},\mathbb{F}_{q^k})$$
$$=\sum\limits_{\stackrel{(u_1,\cdots,u_m)\in
q\setminus(\mathbb{Z}_{(p)}^n/\mathbb{Z}^n)^m}{u_1\otimes v_1+\cdots
+u_m\otimes v_m\equiv \frac{d}{q-1}}}
\sum\limits_{k=1}^{\infty}\frac{-1}{(q^{k|u|}-1)^{mn}}\frac{s^{k|u|}}{k}
\prod_{j=1}^mb_j^{u_j(q^{k|u|}-1)}\prod_{j=1}^m
S_f(T,\omega^{-u_j(q^{k|u|}-1)},\mathbb{F}_{q^{k|u|}})
$$ So
$$C_{f\otimes g,\omega^{-d}}(s,T)
=\exp(\sum\limits_{k=1}^{\infty}\frac{-1}{(q^k-1)^{mn}}
S_f(T,\omega^{-\frac{d}{q-1}(q^k-1)})\frac{s^k}{k})$$
$$=\prod_{\stackrel{(u_1,\cdots,u_m)\in
q\setminus(\mathbb{Z}_{(p)}^n/\mathbb{Z}^n)^m}{u_1\otimes v_1+\cdots
+u_m\otimes v_m\equiv \frac{d}{q-1}}}
\otimes_{j=1}^mC_{f,\omega^{-u_j(q^{|u|}-1)}}(s^{|u|}\prod_{j=1}^m
b_j^{-u_j(q^{|u|}-1)},T,\mathbb{F}_{q^{|u|}}).$$\endproof


\end{document}